\documentclass{fundam}

\publyear{25}
\papernumber{1}
\volume{194}
\issue{1}
\theDOI{10.46298/fi.13734}

\versionForARXIV

\usepackage{amsmath, amscd, amssymb, latexsym, 
}
\usepackage{hyperref}
\usepackage{stackrel}
\usepackage{pdflscape}

\usepackage{tikz}
\usetikzlibrary{matrix,shapes,arrows,positioning,automata}

\def\[#1]{\hbox{$ [\kern -.4em [\, {#1}\,  ]\kern -.4em]\,$}}

\usepackage{float}
\usepackage{xcolor}

\newcommand{\ax}{\displaystyle\mathop{\equiv}_{\mathtt{Ax}}}

\def\[#1]{\hbox{$ [\kern -.4em [\, {#1}\,  ]\kern -.4em]$}}

\newtheorem{theorem}{Theorem}
\newtheorem{lemma}[theorem]{Lemma}
\newtheorem{corollary}[theorem]{Corollary}
\newtheorem{proposition}[theorem]{Proposition}
\newtheorem{definition}[theorem]{Definition}

\newtheorem{remark}[theorem]{Remark}

\newtheorem{example}{Example}

\def\[#1]{\hbox{$ [\kern -.4em [\, {#1}\,  ]\kern -.4em]$}}

\begin{document}

\title{Equational Theory of Ordinals with Addition \\ and
   Left Multiplication by $\omega$}
\author{Christian Choffrut        \\
IRIF, Universit\'{e} Paris Cit\'{e}}

\date{\today}

\maketitle

\runninghead{C.~Choffrut}{Equational theory of ordinals with addition and
   left multiplication by $\omega$}

\begin{abstract}
We  show that the equational theory of the structure
$\langle \omega^{\omega}: (x,y)\mapsto x+y, x\mapsto \omega  x \rangle $
is finitely axiomatizable and give a simple axiom schema when the domain is the set
of transfinite ordinals.
\end{abstract} 

\medskip
\section{Introduction}


Consider   $A_{n}=\{a_{i}\mid 1\leq i \leq n\}$, $n\leq \infty$ and let   
$W_{n}$ denote the set  of  labeled
ordinal words obtained from the empty word and each letter $a$ by closing
under concatenation $\cdot$  ($(u,v)\rightarrow u\cdot v$) and $\omega$-\emph{iteration} ($u^{\omega}= u\cdot u \cdot \cdots$),
example $((ab)^{}c)^{\omega} bc\in W_{3}$. 
Let  $\+W_{n}= \langle W_{n}: 1, \cdot, ^{\omega} \rangle$ denote the resulting algebra.
It satisfies the following infinite (but reasonably simple) set $\Sigma$ of axioms

\begin{align}
\label{eq:left-assoc}
  x \cdot (y \cdot z) &= (x \cdot y) \cdot z \\
\label{eq:left-distributivity}
  (x \cdot y)^{\omega} &= x \cdot (y \cdot x)^{\omega}  \\
  \label{eq:left-domination}
 (x^{p})^{\omega}  &= x^{\omega}\quad p\geq 1\\
   \label{eq:left-unit1}
  x \cdot 1  &= x\\
  \label{eq:left-unit2}
1 \cdot x  &= x
  \end{align}

\bigskip
\noindent
It is known that  $W_{n}$ is isomorphic to the free algebra with $n$ generators in the variety
generated by the equalities
and that for $n>1$,  its equational theory is  axiomatized by the system  $\Sigma$, cf. \cite{BloomChoffrut}. 
The algebra $\+W_{1}$ satisfies identities that do not hold in $\+W_{n}$ with $n\geq2$, for example $xyyx=yxyx$. 
Actually $\+W_{1}$ is isomorphic to the set of ordinals less than $\omega^{\omega}$ 
with the sum and right multiplication by $\omega$ ($x\rightarrow x\omega$).
In \cite{ccseg2020}  it is shown that  for the structure 
$\langle \omega^{\omega}: 0, 1, (x,y)\mapsto x+y, x\mapsto x \omega  \rangle $
(in particular all  constants less than $\omega^{\omega}$ are allowed) equality between two terms is polynomial.
Furthermore the more general algebra of all transfinite ordinals satisfies the same equalities.

Other different  collections of operations on linearly ordered 
labeled words have been studied in the literature. By adding the
$\omega^{\text{opp}}$-iteration ($u^{\omega^{\text{opp}}}= \cdots u\cdot u $) we get the
set $W'_{n}$
 and   the corresponding algebra 
 $\+W'_{n}= \langle W'_{n}: 1, \cdot, ^{\omega}, ^{\omega^*} \rangle$.
 It is  shown that 
for all $n\geq 1$  $\+W'_{n}$ is isomorphic to the free $n$-generated algebra 
in the variety defined by a system of equations which is an axiomatization
of the equational theory
 \cite[Theorem 3.18]{BE2004}.
In \cite{BE2005} a further operation is added, the \emph{shuffle} operator, 
resulting in the set $W''_{n}$ and the corresponding algebra  
$\+W''_{n}= \langle W'_{n}: 1, \cdot, ^{\omega}, ^{\omega^*}, \eta \rangle$. Here again
for all $n\geq 1$  $\+W''_{n}$ is isomorphic to the free $n$-generated algebra 
in the variety defined by a system of equations which is an axiomatization
of the equational theory

Summarising the situation if Eq$(\+A)$ denotes the system of equations satisfied by an algebra $\+A$,
we have Eq$(\+W'_{n})=$Eq$(\+W'_{1})$ and Eq$(\+W''_{n})=$Eq$(\+W''_{1})$ because 
$\+W'_{n}$ and $\+W''_{n}$ are embeddable in $\+W'_{1}$ and $\+W''_{1}$ respectively.
Yet  Eq$(\+W_{n})$ is strictly included in $\text{Eq}(\+W_{1})$ if $n>1$ and an axiom 
schema is  missing for $\+W_{1}$. A semantic approach was considered in \cite{ccseg2020}
Indeed, the right and left handsides of an identity containing $n$ variables
 may be viewed as  mappings of $(\omega^{\omega})^n$ 
into $\omega^{\omega}$. It is proved that  an identity holds in $\+W_{1}$ if and only if the two functions associated with the two handsides  
 coincide over some so-called ``test sets'' such as all the powers $\omega^{n}$ with $n<\omega$ 
along with $0$, for example. Refining this result allows one to exhibit a polynomial time algorithm which determines whether or not two expressions with the  same set of variables define the same mapping, thus hold in $\+W_{1}$. 
It is also proved that over $\omega^{\omega}$ and the transfinite ordinals  the  equational theories are the same.

The algebra investigated here differs from
$\+W_{1}$ in that it considers the left (and not the right) multiplication by $\omega$,
i.e., we consider the signature $\langle (x,y)\mapsto x+y, x\mapsto \omega  x \rangle $. We 
show that it is finitely axiomatizable in $\omega^{\omega}$
and that the more general algebra over the same signature 
but with universe the
transfinite ordinals satisfies a different set of equations of which we give a simple axiom schema.

\section{Ordinals}
We recall the elementary 
properties needed to understand this paper and refer 
 to 
 the numerous standard handbooks such as 
\cite{Rosenstein,Sierpinski} for a more thorough exposition of the theory.
In particular each nonzero ordinal $\alpha$ is uniquely represented by its so-called Cantor normal form.  
$$
\omega^{\alpha_{n}} a_{n} + \cdots + \omega^{\alpha_{0}} a_{0} 
$$
$0<a_{0}, \ldots a_{n}<\omega$ and $\alpha_{n}> \cdots > \alpha_{0}$ 
(a strictly decreasing sequence of ordinals).
The \emph{degree} of $\alpha$ denoted $\partial(\alpha)$ is $\alpha_{n}$,  its
\emph{valuation} denoted $\nu(\alpha)$ is $\alpha_{0}$ and its \emph{length} $|\alpha|$
is the sum  $a_{0} + \cdots +  a_{n}$. The length of $0$ is $0$.

\eject

The sum  $\alpha+\beta$ of two nonnull ordinals $\alpha$ and $\beta$
of   Cantor normal forms
$\alpha=\sum_{i=p}^{i=0}\omega^{\alpha_i}   a_i$ and 
$\beta=\sum_{j=q}^{j=0}\omega^{\beta_j}   b_j$
is the ordinal with the following Cantor normal form 
\begin{align*}
\left(\sum_{i=p}^{i=\ell+1}\omega^{\alpha_i}   a_i\right)
+\omega^{\beta_1}(a_\ell+b_{q}) 
+ \left(\sum_{j=q-1}^{j=0}\omega^{\beta_j}   b_j\right)
& \text{ \ if $\beta_{q}=\alpha_\ell$ for some $\ell$}
\\
\left(\sum_{i=p}^{i=\ell}\omega^{\alpha_i}   a_i\right)
+ \left(\sum_{j=q}^{j=0}\omega^{\beta_j}   b_j\right)
& \text{\ \begin{tabular}{l}
              if $\beta_{q} < \alpha_{\ell }$ for some $\ell$ \\
              and either $\ell=0$ \\
              or $\alpha_{\ell-1}<\beta_{q}$
              \end{tabular}}
\\
\left(\sum_{j=q}^{j=0}\omega^{\beta_j}   b_j\right)
& \text{\ if $\beta_{q} > \alpha_{p }$}
\end{align*}
In particular, 
\begin{align}
\label{eq:degre-valuation-of-sum}
\partial(\alpha+\beta) &= \max(\partial(\alpha), \partial(\beta)) \\
\nonumber
\nu(\alpha+\beta) &=\text{ if $\beta>0$ then $\nu(\beta)$ else $\nu(\alpha)$}
\end{align}

Addition is associative but not commutative:  $\omega = 1+\omega \not= \omega +1$.

If $\alpha+\beta+\gamma=\alpha+\gamma$ we say that $\beta$ \emph{does not count} 
which occurs exactly when $\partial(\beta)<\partial(\gamma)$.

\begin{remark}
\label{re:left-cancelation}
$ \partial(\beta_{1}), \partial(\beta_{2}) < \nu(\alpha_{1})=\nu(\alpha_{2})$ and $ \alpha_{1} + \beta_{1}= \alpha_{2} + \beta_{2}$
implies $ \beta_{1}= \beta_{2}$

\end{remark}

We recall that the left multiplication by $\omega$ is distributive over the sum

\begin{equation}
\label{eq:right-omega}
\omega(\omega^{\alpha_{1}} a_{1} + \cdots + \omega^{\alpha_{n}} a_{n})  =
\omega^{1+ \alpha_{1}} a_{1} + \cdots + \omega^{1+ \alpha_{n}} a_{n} 
\end{equation}
and that $\omega 0=0$.

Every $\alpha\geq \omega^{\omega}$ has a unique decomposition $\alpha = \alpha_{1} + \alpha_{2}$
with $\nu(\alpha_{1})\geq \omega^{\omega}$ and $\partial(\alpha_{2})< \omega^{\omega}$. In this case
because of expression \eqref{eq:right-omega} we have $\omega \alpha_{1} =\alpha_{1}$ and thus
\begin{equation}
\label{eq:omega-omega-omega}
\omega  (\alpha_{1} + \alpha_{2}) = \alpha_{1} + \omega \alpha_{2}
\end{equation}

 \section{Equational axiomatization for $\langle \omega^{\omega}: (x,y)\mapsto x+y, x\mapsto \omega x \rangle $}
\label{sec:S}

We consider the signature with one binary operation  $+$ and one unary 
operation $\omega $ which we interpret  in the  structures 
 $\+S=\langle \omega^{\omega}: (x,y)\mapsto x+y, x\mapsto \omega x \rangle $
and  $\+O=\langle \text{Ord}: (x,y)\mapsto x+y, x\mapsto \omega x \rangle $ whose universes are
respectively  
 the set of all ordinals less
than $\omega^{\omega}$ and the set of all transfinite ordinals.

Let $X$ be a fixed infinite countable set of elements called variables.
The family of \emph{terms} is inductively defined: 
variables  and $0$ are terms, 
if $E$ and $F$ are terms then $E+F$ is a term and 
if $E$ is a  term then $\omega E$ is a term. 
We will avoid unnecessary parentheses by adopting the usual conventions.
We write  $\Sigma_{i=1}^{i=k}E_i$ in place of  $((\ldots(E_1+E_2)+\cdots) +E_{k-1})+E_k$ 
and $\omega^{k} E$ in place of 
$(\omega (\ldots((\omega E\ldots))))$.
In case $k=0$ we let $\omega^0 E$ be $E$.

We write $E(x_{1}, \ldots, x_{n})$ when we need to make explicit the variables on which the term
is inductively constructed. When working in $\omega^{\omega}$
 we will also consider $E(x_{1}, \ldots, x_{n})$ as a function from 
$(\omega^{\omega})^{n}$ into $\omega^{\omega}$  in the natural way. 
Two terms $E$ and $F$ over the same set of variables are \emph{equivalent} 
and we write  $E\equiv F$ 
if they are equal as functions  from 
$(\omega^{\omega})^{n}$  into $\omega^{\omega}$. 
If $\phi$ is an assignment of the variables
in $\omega^{\omega}$   we write $\phi(E)$ to mean $E(\phi(x_{1}), \ldots, \phi(x_{n}))$. 
The notion of being equivalent in $\+O$ is defined similarly.

An \emph{identity} is a pair $(E(x_{1}, \ldots, x_{n}),F(x_{1}, \ldots, x_{n}))$ 
which we write $$E(x_{1}, \ldots, x_{n})=F(x_{1}, \ldots, x_{n})$$ in order to comply with the tradition,
not to be confused with the equality as 
well-formed sequences of symbols over the vocabulary consisting of the variables and $0,+,\omega$.
An identity is \emph{satisfied} in the structure 
if $E$ and $F$ define the same functions.

A set  $\Sigma$ of identities called \emph{axioms} is an \emph{axiomatization} of the equational theory of
$\+S$ if all pairs of equivalent terms $E,F$ are derivable from the axioms by the rules
of equational logic and we write $E\ax F$.

\begin{definition}

 We consider the following set  $\Sigma$ of identities.
 \begin{align}
\label{eq:left-assoc}
  x + (y + z) &= (x + y) + z \\
\label{eq:left-distributivity}
  \omega(x + y) &= \omega x + \omega y  \\
  \label{eq:left-domination}
  x +y + \omega x  &= y + \omega  x\\
  \label{eq:guarded-commutation}
  x +y + z + x+ t+ y  &= y + x + z +x + t+  y\\
  \label{eq:left-unit1}
  x + 0  &= x\\
  \label{eq:left-unit2}
  0 +x  &= x\\
    \label{eq:zero-omega}
  \omega 0&= 0
  \end{align}

\end{definition}

  Let us verify that identity \eqref{eq:guarded-commutation} holds in  $\+S$. 
  Set $\mu= \max (\partial x, \partial y, \partial z, \partial t)$. If  $\partial x, \partial y <\mu$ then 
  both handsides reduce to $ z + x+ t+ y $ because the leftmost occurrences of $x$ and $y$ do not count. 
  If $\partial x =\mu$  and $ \partial y <\mu$ both handsides reduce to 
  $  x  + z + x+ t+ y$ and  if $\partial x <\mu$  and $ \partial y =\mu$ both handsides reduce to 
  $  y  + z + x+ t+ y$. It remains the case where $\partial x = \partial y =\partial (z + x+ t+ y) $ say 
  $x=\omega^{\mu}a + \alpha$, 
  $y=\omega^{\mu}b + \beta$ and $z + x+ t+ y = \omega^{\mu}c + \gamma$
  with  $\partial \alpha,  \partial \beta, \partial \gamma<\mu$. Then both handsides are equal to 
  $\omega^{\mu}(a+b+c )+ \gamma$. 

  The other identities are consequences of the definition 
  of the sum.

Observe the consequence of \eqref{eq:left-domination}
\begin{align}
  \label{eq:c-domination}
 \omega^{p} x +y + \omega^{q}x  &= y + \omega^{q}  x \quad p<q
  \end{align}
Indeed, 
$$
 \omega x^{p} +y + \omega^{q}x =  \omega x^{p} +y + \omega^{p+1}x \cdots + \omega^{q}x
 = y + \omega^{p+1}x \cdots + \omega^{q}x = y + \omega^{q}x 
 $$
\subsection{Elementary properties}

\begin{definition}
  \label{de:x-monomial}
   Given a variable $x$, an $x$-\emph{monomial} 
 is a term of the form $\omega^{e}x$
 where $e$ is a nonnegative integer and $\omega^{e}$ its \emph{coefficient}.
We simply speak of monomial when there is no need to specify which variable.
\end{definition}

\begin{lemma}
  \label{le:normalization}
Via the axioms \eqref{eq:left-distributivity}
and
  \eqref{eq:left-domination}, every term in $\+S$
is equivalent to a term of the form  
\begin{equation}
  \label{eq:normalization}
\sum^{1}_{i=n} \omega^{e_{i}} x_{i}, \quad e_{i}\geq 0
\end{equation}
 with $i>j$ and $x_{i}=x_{j}$ implies $e_{i}\geq e_{j}$. 
   \end{lemma}
  
\begin{proof}
 Trivial by induction and identity \eqref{eq:c-domination}.
 \end{proof}

A term as in \eqref{eq:normalization} is said to be 
\emph{flat}.
Since a flat term can be considered as a sequence of monomials,  the slightly 
abusive expression
``suffix of a term'' has to be understood in the natural way.

\begin{example}
The flat term corresponding to $E=\omega(x+ \omega(\omega x +y + \omega x)+y) +x+x$
is  $$\omega^{3}x+ \omega^{2} y + \omega y + \omega^{3}x + \omega y +x+x$$
which can be identified to the sequence $\omega^{3}x, \omega^{2} y , \omega y , \omega^{3}x, \omega y, x,x$.
 \end{example}

\begin{lemma}
\label{le:invariant-of-subsums}
If two flat terms 
are equivalent, for every $x\in X$ their subsums consisting of the subsequence of their $x$ monomials
 are equal.
   \end{lemma}

\begin{proof}
Put $y=0$ for all $y\not=x$. It suffices to observe that two flat terms with a unique variable $x$ are equivalent if and only if they are equal
since for $x=1$ they evaluate in the same ordinal $1$.
 \end{proof}
 
A \emph{decomposition} $E= E_{k} + \cdots +E_{j} +\cdots +  E_{1}$ of a nonzero flat term $E$
as in   \eqref{eq:normalization} is defined by a subsequence 
$n= i_{k} >i_{k-1} \ldots > i_{1} >i_0=0$ and 
by  grouping successive monomials
\begin{equation}
  \label{eq:grouping}
E_{j} = \sum^{i_{j-1}+1}_{i=i_{j}} \omega^{e_{i}} x_{i}
\end{equation}
We will use two types of decomposition in Lemma  \ref{le:right-most-monomial} and in Theorem 
 \ref{th:decomposition}.

 \begin{lemma}
 \label{le:right-most-monomial}
 Consider  the decomposition of a flat term
 $$
 E=F+\omega^{e}x + G 
 $$
with $e>0$ where $G$ contains no occurrence of $\omega^{e}x$. 
Let $F'+\omega^{e}x + G'$ be the decomposition of an equivalent flat term 
where $G'$ contains no occurrence of $\omega^{e}x$. Then 
$G$   and $G'$ contain the same multiset of occurrences
of $y$-monomials for all  variables $y$.
\end{lemma}
\begin{proof}  Observe that the statement makes sense because by 
Lemma \ref{le:invariant-of-subsums}
if  $E$ and $F$ are equivalent they have the same different monomials. 
It suffices to prove that the set of $y$-monomials in $G$ 
with maximal coefficient 
is an invariant of the equivalence class of $E$.
The statement is clearly true for $y=x$ because we are dealing with flat terms, so we assume $y\not=x$.
Without loss of generality we may assume that the term contains no $z$-monomials
for $z$ different from $x$ and $y$. If there is no occurrence of an $y$-monomial we are done.
Let $\omega^{f}$, $f< e$   be the greatest coefficient of an $x$-monomial
in $G$ if such  an occurrence exist 
and $\omega^{g}$ the greatest coefficient of a $y$-monomial in $G$.
  Set $H= F+\omega^{e}x $ 
and let $a,b$ be two integers such that
$a+f=b+g$. Then $\nu(H(\omega^{a}, \omega^{b}))=e+a  $
and $G(\omega^{a}, \omega^{b}) =\omega^{a+f}\cdot n + \alpha$ where
$\partial(\alpha)< a+f$ and
$n$ is the number of occurrences of $\omega^{f}x$ plus the number of occurrences of  $\omega^{g}y$.
If $G$ contains no $x$-monomial, then 
$\nu(H(\omega^{g}, 1))=e+g  $ and 
$G(\omega^{g}, 1)=\omega^{g}\cdot n + \alpha$ where
$\partial(\alpha)<g$ and $n$ is the number of occurrences of $\omega^{g}y$.
\end{proof}

\begin{definition}
\label{de:new-monomial}
With the notations of   \eqref{eq:normalization} the 
 \emph{new monomial  decomposition} (NMD) of $E$ is the sum
$E=E_{n}+E_{n-1}+ \cdots + E_{1}$
 where $E_{i}=E'_{i}+ \omega^{e_{i}} x_{i}$ for some $E'_{i}$ such that  $\omega^{e_{i}}x_{i}$ does not occur
in $E_{i-1}+\cdots + E_{1}$ (by convention $E_{0}=0$). 
\end{definition}

Observe that $n$ is the number of different monomials in $E$. The idea is to record 
the moment when  a new monomial appears  in a 
scan from right to left.

\begin{lemma}
 \label{le:new-monomial-decompostion} 
Let $E$ and $F$ be two flat terms and their  NMD
\begin{equation}
\label{eq:new-monomial-decomposition}
\begin{array}{l}
E= E'_{n} +\omega^{e_{n}}x_{n}+\cdots + E'_{1} +\omega^{e_{1}}x_{1}\\
F= F'_{m} +\omega^{f_{m}}y_{m}+\cdots + F'_{1} +\omega^{f_{1}}y_{1}
\end{array}
\end{equation}
If  $E\equiv F$ then  $n=m$, $\omega^{e_{i}}x_{i}=\omega^{f_{i}}y_{i}$
for $i=1, \ldots, n$ and 
for $i=1, \ldots n$ $E'_{i}$ and  $F'_{i}$
differ by a permutation of their monomials.
 
 \end{lemma}

 \begin{proof}
Clearly $n=m$ because  $E$ and $F$ have the same different monomials by 
Lemma \ref{le:invariant-of-subsums}.
The last claim is a consequence of   Lemma   \ref{le:right-most-monomial}.
 \end{proof}
    \begin{theorem}
  \label{th:decomposition}
$\Sigma$ is an axiomatization of $\+S$.
\end{theorem}

\begin{proof}
Consider the  new monomial decompositions
\eqref{eq:new-monomial-decomposition}. 
By Definition  \ref{de:new-monomial}  for all successive monomials in $E'_{i}$ (resp. in  $F'_{i}$)
there exists an occurrence of these monomials in $\omega^{e_{i}}x_{i} + \cdots + E'_{1} +\omega^{e_{1}}x_{1}$
resp. in  $\omega^{e_{i}}x_{i} + \cdots + F'_{1} +\omega^{e_{1}}x_{1}$. Since each permutation is  a product of transpositions,
 we get  $E\ax F$ by repetitive applications of axiom \eqref{eq:guarded-commutation}.
\end{proof}

 \section{Equational axiomatization for 
 $\langle \+O: (x,y)\mapsto x+y, x\mapsto \omega x \rangle$}

Property \eqref{eq:guarded-commutation} is valid in $\+O$,
by interpreting in the proof of Section \ref{sec:S}, the exponent $\mu$ 
as an element of  $\+O$.
 Beyond $\omega^{\omega}$  axiom \eqref{eq:left-domination}  no longer holds. Indeed
 since $\omega \omega^{\omega} = \omega^{1+ \omega} = \omega^{\omega} $.
 \begin{equation}
\label{eq:beyond}
 \omega^{\omega}  \cdot 2 =
  \omega^{\omega} + \omega \omega^{\omega}\not=  \omega \omega^{\omega}= \omega^{\omega} 
\end{equation}
We consider the system $\Sigma'$ consisting of the axioms of $\Sigma$
except axiom \eqref{eq:left-domination} along with the following new axioms.


%
%
\begin{equation}
\label{eq:domination1}
 \omega^{p} x +y + \omega^{q}x  = x + y + \omega^{q}  x \quad 0< p<q
\end{equation}
\begin{equation}
\label{eq:hidden-old-variable1}
  x + \omega^{r}y + t + \omega^{p}x  + u + \omega^{q}y =  \omega^{r}y + x +  t + \omega^{p}x  + u + \omega^{q}y    \quad p>0  \text{ or } q=r ,
\end{equation}
\begin{equation}
\label{eq:hidden-old-variable2}
  x + \omega^{r}y + t  + \omega^{q}y  + u + \omega^{p}x  =  \omega^{r}y + x +  t +  \omega^{q}y + u + \omega^{p}x      \quad p>0  \text{ or } q=r  \end{equation}

\begin{lemma}
Identities \eqref{eq:domination1}, \eqref{eq:hidden-old-variable1} 
 hold in the structure $\langle \+O: (x, y)\mapsto x+y, x\mapsto \omega x \rangle$. 
\end{lemma}

\begin{proof}
 Indeed, \eqref{eq:domination1} holds in $\omega^{\omega}$ since by 
 identity    \eqref{eq:left-domination} both hand sides reduce to $y+\omega^{q}x$. If  $x\geq \omega^{\omega}$ 
we decompose 
 $x= x_{1} + x_{2} $ with 
 $\nu x_{1}\geq \omega$ and $\partial x_{2}<\omega$ as in  expression  \eqref{eq:omega-omega-omega}.
 Because $\omega^{p} \cdot x_{1}= x_{1}$ we get 
 $\omega^{p} \cdot x = x_{1} + \omega^{p} \cdot x_{2}$
and thus
 $$
 \begin{array}{lll}
  &\omega^{p} x  +y  + \omega^{q}x  \\
= &  x_{1} + \omega^{p} \cdot x_{2}  + y + 
	 \omega^{q} \cdot x & \text{(decomposition of } x)\\
  = &
 x_{1} +  y + 
 \omega^{q} \cdot x
	 & (\partial \omega^{p} x_{2} <  \partial  \omega^{q}   x)\\
   = &
 x_{1} +  x_{2}  + y + 
 \omega^{q} \cdot x
	 & (\partial x_{2} <  \partial  \omega^{q}   x)  \\
=&  x + y + 
	 \omega^{q} \cdot x & \text{(recomposition of } x)
  \end{array}
   $$

 Concerning \eqref{eq:hidden-old-variable1} and \eqref{eq:hidden-old-variable2},
  if $\partial (x), \partial(y)\geq \omega$ in both cases the two hand sides reduce 
$$  x +  y + t +  x  + u +  y =   y + x +  t +  x  + u + y$$ 
which holds because of identity \eqref{eq:guarded-commutation}.
%
%
  Consider   \eqref{eq:hidden-old-variable1}. If $p=0$, thus $q=r$,  then this is  axiom \eqref{eq:guarded-commutation}.
If $p>0$ and $x<\omega^{\omega}$ 
 the leftmost $x$ does not count 
and if $x\geq \omega^{\omega}$ then $\omega^{r} y$ does not count. Identity \eqref{eq:hidden-old-variable2} is proved similarly.
\end{proof} 
Because of identity \eqref{eq:domination1}  every term is equivalent to a term of the form
\begin{equation}
  \label{eq:normalization1}
E= \sum^{1}_{k=m} \omega^{e_{k}} x_{k}, \quad e_{k}\geq 0
\end{equation}
 where $i>j$,  $x_{i}=x_{j}$ and $e_{i}, e_{j}>0$ implies $e_{i}\geq e_{j}$
 which we  call \emph{pseudo flat}. An occurrence of a variable $x$ (i.e., a monomial 
 with coefficient $1$)  to the left of some monomial
 $\omega^{e}x $ with $e>0$ is called  \emph{hidden}.
 \begin{example}
 \label{ex:normalization}
In 
$$
z + \omega x + \underline{x}+ \underline{y} +\omega y + \omega x + y + \omega x + y + x
$$
the unique hidden occurrences are the leftmost underlined occurrences of $x$ and $y$.
\end{example}
The definition \ref{de:new-monomial} of new monomial  decomposition 
extends naturally to the present structure. 
\begin{example}
(Example \ref{ex:normalization} continued) 
With the notations of 
definition \ref{de:new-monomial}  we have 
$$
(z) + (\omega x + x+ y +\omega y) + (\omega x + y + \omega x) + (y) + (x)
$$
Indeed, when reading the expression from right to left we define five groups of terms and start a group whenever 
we find a monomial that never occurred before (the rightmost monomial of each group). For example the third group (from the right) starts
	with the new monomial $\omega x$
and contains no new monomial; it ends just  before a new monomial ($\omega y$) appears.
Observe that the hidden occurrences are not deleted.
\end{example}
The following is obtained by simple adaptation of Lemma  \ref{le:new-monomial-decompostion}.
 \begin{lemma}
 \label{le:new-monomial-ordinal}
 Consider the new monomial decomposition of  two equivalent  
pseudo-flat terms $E$ and $F$.
 Then
 $$
 \begin{array}{l}
 E= E'_{n} +\omega^{e_{n}}x_{n} + \cdots + E'_{1} +\omega^{e_{1}}x_{1}\\
 F= F'_{n} +\omega^{e_{n}}x_{n} + \cdots + F'_{1} +\omega^{e_{1}}x_{1}\\
 \end{array}
  $$
 where for all $i=1, \ldots, n$   $E'_{i}$ and $F'_{i}$ contain the same occurrences
 of nonhidden monomials.
 \end{lemma}
 Now instead of splitting accordingly to a ``new monomial'' we split according to a ``new variable'' (have in mind a scan from right to left). 
 \begin{definition}
\label{de:new-variable}
With the notations of   \eqref{eq:normalization1} the 
 \emph{new variable  decomposition} (NVD) of a pseudo-flat $E$ is the sum
$E=E_{n}+E_{n-1}+ \cdots + E_{1}$
 where for each $i=1, \ldots, n$  $E_{i}=E'_{i}+ \omega^{e_{i}} x_{i}$ such that  no $x_{i}$-monomial appears 
in $E_{i-1}+\cdots + E_{1}$ with the convention $E_{0}=0$. 
\end{definition}
No hidden occurrence can be the right-most monomial of some $E_{i}$
and  that the number of subterms in the decomposition equals the number of variables in the term.

\begin{example} 
(Example \ref{ex:normalization} continued). The new variable decomposition is a sum of three subterms.
$$(z) + (\omega x + x+ y +\omega y + \omega x + y + \omega x + y) + (x)$$
\end{example}
Observe the difference with the new monomial decomposition.
Here the expression is decomposed  
in three groups. From right to left: 
the $x$-monomials, then the $y$-monomials with possibly some $x$-monomials and
finally the  $z$-monomials  and possibly some $y$- or $x$-monomials.
here 

\begin{lemma}
\label{le:hidden-occurrecnes-in-ord}
Let $E(x_{1}, \ldots, x_{n})$ and $F(x_{1}, \ldots, x_{n})$ be two equivalent pseudo-flat terms  and 
consider their new variable  decompositions $E=E_{n} + E_{n-1} + \cdots + E_{1}$
and $F=F_{m} + F_{n-1} + \cdots + F_{1}$. Then
  $n=m$
and for all $i\leq j$  the number of hidden occurrences  of each $x_{i}$ in $E_{j}$ is equal to the number of hidden occurrences 
 $x_{i}$ in $F_{j}$. 
\end{lemma}
\begin{proof} Equality $n=m$ is obvious because it is the number of variables of the two equivalent terms $E$ and $F$.
By possibly renaming the variables we may assume that for $i=1, \ldots, n$  the rightmost monomial in $E_{i}$
is $\omega^{e_{i}} x_{i}$ for some $e_{i}\geq 0$. By definition a hidden occurrence of $x_{i}$ can only belong to the 
set of monomials in $E_{n}+ \cdots + E_{i}$. We set for all $i\leq j\leq n+1$
\begin{itemize}
\item $H_{i,j}$ the number of hidden occurrences of  $x_{i}$ in $E_{j-1}+ \cdots + E_{i}$

\item $N_{i,j}$ the number  occurrences of nonhidden  $x_{i}$-monomials in $E_{j-1}+ \cdots + E_{i}$
\end{itemize}
 Consider the assignment $\phi(x_{i})=1$,
$\phi(x_{j})=\omega^{a}$ and $\phi(x_{\ell})=0$ for $\ell\not=i,j$.  For  $a<\omega $ greater than the exponents of
the coefficients of 
all $x_{i}$-monomials 
we have
$
\phi(E) = \alpha + \beta 
$
where $\nu(\alpha)\geq a> \partial(\beta) $. Then $N_{i,j}=|\beta|$. Now consider the assignment
$\phi(x_{i})=\omega^{\omega}$, $\phi(x_{j})=\omega^{\omega +1}$ and $\phi(x_{\ell})=0$ for $\ell\not=i,j$. 
Then we have $
\phi(E) = \omega^{\omega +1} c +  \omega^{\omega}d, c,d<\omega $ and $H_{i,j}=d- N_{i,j}$. Then the number of hidden occurrences of $x_{i}$ in $E_{j}$
is equal to $H_{i,j+1} - H_{i,j}$
\end{proof}
\begin{theorem}
The system of identities $\Sigma'$ is an axiomatization of the equational theory 
of  the structure $\langle \text{\em Ord}: (x,y)\mapsto x+y, x\mapsto \omega x \rangle $
\end{theorem}
\begin{proof}
Let $E$ and $F$ be two equivalent pseudo-flat terms and let 
 $E=E_{n} + E_{n-1} + \cdots + E_{1}$
and $F=F_{n} + F_{n-1} + \cdots + F_{1}$ be their NVD.
We show 
that we have $E\ax \widehat{E}$ where $\widehat{E}=F_{n} + E_{n-1} + \cdots + E_{1}$.
By Lemma \ref{le:hidden-occurrecnes-in-ord} $E_{n}$ and $F_{n}$
have the same number of hidden occurrences for all variables. We now consider their
NMD. By Lemma \ref{le:new-monomial-ordinal}
we have
 $$
 \begin{array}{l}
 E_{n}= E'_{n,s} +\omega^{e_{n,s}}x_{n,s} + \cdots + E'_{n,1} +\omega^{e_{n,1}}x_{n,1}\\
 F_{n}= F'_{n,s} +\omega^{e_{n,s}}x_{n,s} + \cdots + F'_{n,1} +\omega^{e_{n,1}}x_{n,1}\\
 \end{array}
  $$
where for $j=s, \ldots , 1$,  $E'_{n,j}$ and $F'_{n,j}$ have the same set of occurrences 
of nonhidden monomials. By equations \eqref{eq:hidden-old-variable1},  \eqref{eq:hidden-old-variable2} 
(with $r=q$) and  \eqref{eq:guarded-commutation}, 
all occurrences of nonhidden monomials in $E'_{n,j}$ commute pairwise and with the hidden monomials.
So we have
with $R= E_{n-1} + \cdots + E_{1}$
$$
\begin{array}{rl}
& E'_{n,s} +\omega^{e_{n,s}}x_{n,s} + \cdots + E'_{n,1} +\omega^{e_{n,1}}x_{n,1}+ R\\
\ax  & E''_{n,s} +\omega^{e_{n,s}}x_{n,s} + \cdots + E''_{n,1} +\omega^{e_{n,1}}x_{n,1}+ R\\
\end{array}
$$
where for $j=s, \ldots, 1$, 
$E''_{i,j}$ and  $F'_{i,j}$ have the same sequence of nonhidden monomials. They may only differ 
in the number and positions of the hidden monomials. 
Since all hidden monomials commute by \eqref{eq:guarded-commutation} with all monomials 
in $E''_{n,s} +\omega^{e_{n,s}}x_{n,s} + \cdots + E''_{n,1} +\omega^{e_{n,1}}x_{n,1}+ R$ 
and since the number of hidden monomials in $F_{n}$ and $E''_{n,s} +\omega^{e_{n,s}}x_{n,s} + \cdots + E''_{n,1} +\omega^{e_{n,1}}x_{n,1}$ are the same for all variables,
we have
$$
\begin{array}{rl}
   F_{n}  \ax E''_{n,s} +\omega^{e_{n,s}}x_{n,s} + \cdots + E''_{n,1} +\omega^{e_{n,1}}x_{n,1}+ R
\end{array}
$$
Now $E\ax \widehat{E}$ implies $E\equiv \widehat{E}$ thus $\widehat{E} \equiv F$ and 
finally by cancelation  $ E_{n-1} + \cdots + E_{1}\equiv  F_{n-1} + \cdots + F_{1}$ which allows us
to conclude by induction.
\end{proof}

\section{Complexity}

We show that  the complexity of determining the equivalence of two terms $E$ and $F$ is linear in the size of the expressions
as element of the algebras $\+S$ and $\text{Ord}$ respectively.

The flattening of a term can be achieved as follows. Construct in linear time the syntactic tree whose nodes are of arity $1$ 
for the left multiplication by $\omega$ and of arity $2$ for the addition. The leaves of the tree are labeled by the different variables occurring 
in the term. Then a depth-first search assigns to each variable its coefficients $\omega^{i}$. 

A flat expression is a sum of monomials. 
Thus we may consider the finite set  $A$ whose elements $a_{ij}$ are in one-to-one correspondence with the different
monomials $\omega^{e_{j}} x_{i}$, $1\leq i\leq n, 1\leq j \leq m_{i}$ occurring in the term. 
Two terms can only be equivalent if they have the same multiset of monomials.  We define $A_{i}=\{a_{i,j}: 1\leq j\leq m_{i}\}$
for all $1\leq i\leq n$ and thus $A=\bigcup^{n}_{i=1} A_{i}$. A flat expression can 
be identified with a sequence $u\in A^*$ satisfying the conditions
\begin{itemize}
\item for all $i$ there exists $1\leq j\leq  m_{i}$ and $v,w$ such that $u=va_{ij}w$ (by the definition of $A$)
\item for all factorizations $u=va_{ij} w a_{ik}z$ it holds $j\geq k$. (see  condition \eqref{eq:normalization})
\end{itemize}
For all $v\in A^*$ we set  $\text{C}(v)= \{a_{ij}\in A: \exists w,z\ (v=wa_{ij}z) \}$. 
Henceforth by flat expression we mean a sequence on $A$ subject to the above conditions. 
A flat expression has a 
unique factorization $u=u_{n}u_{n-1} \cdots u_{1}$, which we call  \emph{subalphabet factorization} (see the NMD) 
 such that $u_{1}$ is the longest word such that $|C(u_{1}|=1$  and  for all $i>1, u_{i}$ is the longest  word such that $|C(u_{i}\cdots u_{1})\setminus C(u_{i-1}\cdots u_{1})|=1 $

The equivalence of two flat terms in $\+S$ can be done by comparing in a single pass
their subalphabet factorizations 
$u=u_{n}u_{n-1} \cdots u_{1}$  and $v_{m}v_{m-1} \cdots v_{1}$
and by checking whether or  not $C(u_{i})=C(v_{i})$ for all $i\leq 1$. 
The equivalence of two pseudo-flat terms in  $\text{Ord}$
can be done by.  additionally counting the hidden variables in the same pass.

The following shows that $E$ and $F$ are equivalent if and only if their restrictions to any pair of variables are equivalent
(the restriction to $\{x,y\}$ consists in deleting all $z$-monomials for $z$ different from $x$ and $y$). 

\begin{proposition}%
Two expressions $E$ and $F$ are equivalent if and only if their restrictions on  any two pair of variables
are equivalent.
\end{proposition}

\begin{proof}%
We  identify  two expressions with two sequences  $u,u'\in A^*$ as above.
Consider their  subalphabet factorizations $u=u_{n}u_{n-1} \cdots u_{1}$
and $u'=u'_{n}u'_{n-1} \cdots u'_{1}$ (they have the same length).
 For every 
$1\leq i,j\leq n$  let $\pi_{i,j}$ be the projection of $A^*$ onto $(A_{i}\cup A_{j})^*$.
The  rightmost letter of  $u_{i}$ and $u'_{i}$ are the same. 
This is due to the fact that for all pairs $(i,j)$ we have
$\pi_{i,j}(u)=\pi_{i,j}(u')$ showing that all rightmost occurrences of $a\in A$ in $u$ and $u'$
we have $u=vaw$ and $u'=v'aw'$ with $C(w)=C(w')$. 

Now we prove that for all $a\in A_{k}$ 
the number of occurrences of  $a$ in $u_{i}$ depends only on the projections $\pi_{i,\ell}$ 
for $1\leq \ell \leq n$.
Indeed, assume  
the rightmost letters of $u_{i}$ and $u_{i+1}$ are in $A_{r}$ and $A_{s}$ 
respectively. 
The number of occurrences of $a$ in $u_{n} \cdots u_{i}$
is equal to the number of occurrences of $a$ in $\pi_{k,s}(u_{n}u_{n-1} \cdots u_{i})$
and the number of occurrences of $a$ in $u_{n} \cdots u_{i+1}$
is equal to the number of occurrences of $a$ in  $\pi_{k,s}(u_{n}u_{n-1} \cdots u_{i+1})$. 
\end{proof}

\begin{corollary}%
Equality $E=F$ is provable if and only if for all $1\leq i<j\leq n$
equalities $E_{ij}=F_{ij}$ are provable.
\end{corollary}


\begin{thebibliography}{99}
%
\bibitem{BloomChoffrut}
S. L. Bloom and C. Choffrut,
Long words: the theory of concatenation and $\omega$-power, 
{\em Theoret. Comp. Sci.}, 259 (2001), 533--548. 

\bibitem{ccseg2020}
C. Choffrut and S. Grigorieff,
Test sets for equality of terms in the additive structure of ordinals augmented with right multiplication by $\omega$,
{\em Algebra Universalis}, 81, 45, 2020. 



 \bibitem{BE2004}
S. L. Bloom and
                  Z. {\'{E}}sik.
\newblock Axiomatizing omega and omega-op powers of words.
\newblock {\em  RAIRO Theor. Informatics Appl.}, 38, 1, 3--17, 2004.
  
  
 \bibitem{BE2005}
S. L. Bloom and
                  Z. {\'{E}}sik.
\newblock The equational theory of regular words.
\newblock {\em Inf. Comput.},  197, 1-2, 55--89, 2005.


\bibitem{Rosenstein} 
J.G. Rosenstein, 
{\em Linear ordering}, 
Academic Press, New-York, 1982

\bibitem{Sierpinski} 
W. Sierpi{\'n}ski,
{\em Cardinal and Ordinal Numbers}, 
2nd ed., PWN-Polish Scientific Publishers, 
1965.



\end{thebibliography}
\end{document}